\documentclass[11pt]{article}
\usepackage{amsfonts,amsmath,amscd,amssymb,amsthm, graphics}

\newfont{\eufm}{eufm10 scaled\magstep1}
\newcommand{\FRAK}[1]{\mbox{\eufm#1}}

\newcommand{\cO}{\mathcal{O}}
\newcommand{\cC}{\mathcal{C}}
\newcommand{\cI}{\mathcal{I}}

\newcommand{\cD}{\mathcal{D}}
\newcommand{\cB}{\mathcal{B}}
\newcommand{\cV}{\mathcal{V}}

\newcommand{\cP}{\mathcal{P}}

\newcommand{\cH}{\mathcal{H}}

\newcommand{\cW}{\mathcal{W}}

\newcommand{\bbN}{\mathbb{N}}
\newcommand{\bbZ}{\mathbb{Z}}

\newcommand{\bbR}{\mathbb{R}}

\newcommand{\bbF}{\mathbb{F}}

\newcommand{\bbX}{\mathbb{X}}
\newcommand{\bbY}{\mathbb{Y}}

\newtheorem{thm}{Theorem}[section]
\newtheorem{lem}[thm]{Lemma}

\newtheorem{prop}[thm]{Proposition}

\newtheorem{rem}[thm]{Remark}

\newtheorem{al}[thm]{Algorithm}
\begin{document}
\title{On the computation of graded components of Laurent polynomial rings}
\author{Sonia L. Rueda \thanks{partially supported by the Spanish ``
Ministerio de Educaci\'on y Ciencia" under the Project
MTM2005-08690-C02-01 and  by the ``Direcci\'on General de
Universidades de la Consejer\'{\i}a de Educaci\'on de la CAM y la
Universidad de Alcal\'a" under the project CAM-UAH2005/053. }\\
Departamento de Matem\' atica Aplicada, E.T.S. Arquitectura\\
Universidad Polit\' ecnica de Madrid\\
Avda. Juan de Herrera 4, 28040-Madrid, Spain\\
{sonialuisa.rueda@upm.es}}

\maketitle

\begin{abstract}
In this paper, we present several algorithms for dealing with
graded components of Laurent polynomial rings. To be more precise,
let $S$ be the Laurent polynomial ring
$k[x_1,\ldots,x_{r},x_{r+1}^{\pm 1},\ldots ,x_n^{\pm 1}]$, $k$
algebraicaly closed field of characteristic $0$. We define the
multigrading of $S$ by an arbitrary finitely generated abelian
group $A$. We construct a set of fans compatible with the
multigrading and use this fans to compute the graded components of
$S$ using polytopes. We give an algorithm to check whether the
graded components of $S$ are finite dimensional. Regardless of the
dimension, we determine a finite set of generators of each graded
component as a module over the component of homogeneous
polynomials of degree 0.
\end{abstract}

%%%%%%%%%%%%%%%%%%%%%%%%%%%%%%%%%%%%%%%%%%%%%%%%%%%%%%%%%%%%%%%%%%%
\section{Introduction}
%%%%%%%%%%%%%%%%%%%%%%%%%%%%%%%%%%%%%%%%%%%%%%%%%%%%%%%%%%%%%%%%%%%

The formal construction of the homogeneous coordinate ring of a
toric variety, that I.M. Musson in \cite{M2}, D.A. Cox in
\cite{Cox} and others discovered in the early 1990s, takes a fan
$\Delta $ and creates a torus action on an open subset of an
affine space whose quotient is the toric variety of $\Delta $. We
reverse this process in \cite{R}. To be more specific, let $k$ be
an algebraically closed field of characteristic $0$. We take a
torus action on the affine space $X=k^r\times (k^{\times})^s$ and
create various fans whose toric varieties are the quotients under
the action of the torus of an open subset of $X$.  This "reverse
engineering" has appeared in other places, see for instance
Chapter 10 of Miller and Sturmfels book \cite{MS} where they take
an algebraic torus times a finite abelian group acting on $k^n$.

In this paper, we extend the construction of the fan associated to
the action given in \cite{R} to cover also Miller and Sturmfels
situation but we focus on making the process computational. Let
$H$ be an algebraic torus times a finite abelian group, we give an
algorithm to obtain a set of fans associated to the diagonal
action of $H$ on $X$. Using these fans, we will be able to give a
computational description of the graded components of Laurent
polynomial rings.

Let $S$ be the Laurent polynomial ring
$k[x_1,\ldots,x_{r},x_{r+1}^{\pm 1},\ldots ,x_n^{\pm 1}]$, with
$n=r+s$. The action of $H$ on $X$ extends to an action of $H$ on
$S$. As in the polynomial ring case, this action determines a
multigrading of $S$ by the finitely generated abelian group
$A=Hom(H,k^{\times})$. Multigradings of polynomial rings are
treated in detail in \cite{MS}, Chapter 8.  In section
\ref{multigrading}, we extend the definition of multigrading to
Laurent polynomial rings and in section \ref{compFan}, we provide
an algorithm to obtain a set of fans compatible with the
multigrading of $S$ by $A$. This fans allow us to see graded
components of $S$ as polyhedra. Each graded component $S_a$, $a\in
A$ is a module over the ring of differential invariants $\cD
=k[x_1,\ldots,x_{r},x_{r+1}^{\pm 1},\ldots ,x_n^{\pm 1},
\partial_1 ,\ldots ,\partial_n]^H$ where $\partial_i=\partial/\partial
x_i$, $i=1,\ldots ,n$. In \cite{MR}, we studied graded components
of $S$ from this point of view . The connection between group
actions and finite fans allows a computational description of this
family of $\cD$-modules that we present in this paper.

When computing the graded components of $S$, it is important to
decide whether they are finite dimensional as $k$ vector spaces or
not. Let $S_0$ be the $k$ vector space of homogeneous polynomials
of degree $0$. As in the polynomial ring case, either all the
graded components are finite or infinite dimensional as we
conclude in the last section of the paper. In section
~\ref{special}, we give an algorithm to check whether $S_0$ is
finite dimensional.

Given a finite fan $\varDelta$, a shift with nonempty intersection
of the dual cones of the cones of $\varDelta$  gives a polyhedron.
Furthermore, such a polyhedron is a polytope if and only if the
fan is not contained in a half-space, as we prove in section
~\ref{fans}. In section ~\ref{polytopes}, given a fan compatible
with the multigrading of $S$ by $A$, we use polyhedra obtained in
this manner to give a set of generators for each $S_a$ as a module
over $S_0$. If the graded components of $S$ are finite dimensional
we compute a $k$-basis of $S_a$ using the lattice points of a
polytope. We provide a set of linear inequalities determining this
polytopes so that the number of lattice points can be counted with
LattE, \cite{LTY}. When $S_0$ is infinite dimensional,  we
determine a finite set of generators $f_1,\ldots ,f_h$,
corresponding to the Hilbert basis of a cone, such that
$S_0=k[f_1,\ldots ,f_h]$. To determine the infinite dimensional
$S_a$, we compute a finite set of generators of $S_a$ as an
$S_0$-module corresponding to the lattice points of a polytope.

We begin our paper describing in section ~\ref{action} the group
$H$ and the diagonal action of $H$ on $X$. Throughout this paper,
a particular action will be specified by a matrix $L$ whose
columns are the weights of the action. In some special cases, we
can replace the matrix $L$ by another one with the special form
given in section ~\ref{specialmatrix} which gives the same action.
Most of the algorithms given in this paper will take the matrix
$L$ or its special form as an input.

%%%%%%%%%%%%%%%%%%%%%%%%%%%%%%%%%%%%%%%%%%%%%%%%%%%%%%%%%%%%%%%%%%%%
\section{Description of the group and the action}\label{action}
%%%%%%%%%%%%%%%%%%%%%%%%%%%%%%%%%%%%%%%%%%%%%%%%%%%%%%%%%%%%%%%%%%%%

We consider a diagonal action of $H$ on $X=k^r\times
(k^{\times})^s\subseteq k^n$ with $n=r+s$. This is an action that
extends to a diagonal action on $k^n$. Such an action is given by
an embedding of $H$ into the group $T$ of diagonal matrices in
$GL_n(\bbZ)$. Details about this action in the algebraic torus
case are given in \cite{MR}, \S 2.1. The construction developed
there applies to the more general case treated here.

Identify $H$ with $G\times\bbF$ where $G=(k^{\times})^p$ and
$\bbF$ is the product of cyclic groups of orders $d_1,\ldots
,d_t$. We can identify the group of characters of $H$,
$A=Hom(H,k^{\times})$ with
\begin{displaymath}
A=\bbZ^p\oplus \bbZ/d_1\bbZ \oplus\ldots\oplus\bbZ/d_t\bbZ.
\end{displaymath}
We think of $A$ as a space of column vectors with integer entries.

There exist $\eta_1,\ldots ,\eta_n\in A$ such that {\em H acts on
$X$ with weights $\eta_{1},\ldots ,\eta_n$ .} Let $m=p+t$. We
denote by $L$ the $m\times n$ matrix with $i$-th column vector
$\eta_i$, $i=1,\ldots ,n$. We say that {\em H acts on $X$ by the
matrix $L$}. We may and we will assume that $H$ acts faithfully on
$X$. Therefore $L$ has rank $m$.

Let us also denote by $D$ the diagonal matrix with entries
$d_1,\ldots ,d_t$. This $t\times t$ matrix will be used in several
computations throughout this paper.

%%%%%%%%%%%%%%%%%%%%%%%%%%%%%%%%%%%%%%%%%%%%%%%%%%%%%%%%%%%%%%%%%%%%
\section{Finite fans}\label{fans}
%%%%%%%%%%%%%%%%%%%%%%%%%%%%%%%%%%%%%%%%%%%%%%%%%%%%%%%%%%%%%%%%%%%%

In this section we state some results about finite fans that will
be needed in the rest of the paper. As far as possible we follow
the notation of \cite{F} , Chapter 1.

Let $N= \bbZ^l$ be the $l$-dimensional lattice. Let $\varDelta $
be a fan in $N$, which is a collection of strongly convex rational
polyhedral cones in $N_{\bbR}=N\otimes _{\bbZ} \bbR\equiv\bbR^l$.
Through this paper, cones will be strongly convex, rational and
polyhedral except in a few occasions where we will say convex cone
to mean convex rational polyhedral cone.

Let $M=Hom_{\bbZ}(N,\bbZ)$ and $\langle \_\;,\_\rangle : M\times
N\rightarrow \bbZ$ the natural bilinear pairing. For each
$\sigma\in\varDelta$, let
\begin{equation}\label{lambda}
\Lambda_{\sigma}=M\cap\sigma^{\vee}=\{u\in M|\langle u,v\rangle\geq
0\mbox{ for all }v\in\sigma\}.
\end{equation}

We review some terminology introduced in \cite{R}, \S 1.2. Denote
by $\varDelta (1)$ the set of cones of $\varDelta $ with dimension
one. Given $v\in N$, let $\tau_v=\bbR_+ v$ be the ray generated by
$v\in N$. Given $\sigma\in\varDelta$ we define $[\sigma
]=\{i\in\{1,\ldots ,r\} \mid \tau_{v_i}\mbox{ is a face of }\sigma
\}$. Given $u\in M_{\bbR}=M\otimes_{\bbZ}\bbR\equiv\bbR^l$, a
subset of the form $H_{u}=\{ v\in N_{\bbR} \mid \langle
u,v\rangle\geq 0\}$ with $u\neq 0$ is called a half-space in
$N_{\bbR}$. We say that the fan $\varDelta $ is {\em contained in
a half-space} if we can find $0\neq u\in M_{\bbR}$ such that
$\sigma\subseteq H_{u}$ for all $\sigma\in\varDelta$.

%%%%%%%%%%%%%%%%%%%%%%%%%%%%%%%%%%%%%%%%%%%%%%%%%%%%%%%%%%%%%%%%%%

Given a fan $\varDelta$, let $\{v_1,\ldots ,v_r\}$ be a set of
vectors in $N$ generating the one dimensional cones of
$\varDelta$. Fix $\varphi=(\varphi_1,\ldots
,\varphi_n)\in\bbN^r\times\bbZ^s$ and define the sets
\begin{equation}\label{gamma-sigma-m}
\gamma_{\sigma ,\varphi}=\{u\in M_{\bbR} \mid \langle
u,v_{i}\rangle\geq -\varphi_i \mbox{ for all } i\in [\sigma] \}
\end{equation}
for each $\sigma\in\varDelta$. Consider the polyhedron
\begin{equation}\label{polytope}
\cP_{\varphi}=\displaystyle{\cap_{\sigma\in\varDelta}}\;\gamma_{\sigma
,\varphi}= \{u\in M_{\bbR} \mid \langle u,v_i\rangle \geq
-\varphi_i\mbox{ for all }i=1,\ldots ,r\}.
\end{equation}
The purpose of the following two lemmas is to prove that
$\cP_{\varphi}$ is bounded if and only if the fan $\varDelta$ is not
contained in a half-space.

\begin{lem}\label{continuous}
Fix an integer $a\geq 0$ and $v\in\bbZ^l$ and consider the set
$\Omega =\{u\in\bbR^{l} \mid \langle u,v\rangle\geq -a \} $. Given
$u\in S=\{u\in \bbR^{l}||u|=1\}$, let $l_{u}=sup \{\mu\in\bbR^{+}
| \mu u\in \Omega\}$. Then the map $f:S\rightarrow \bbR^{l}$
defined by $f(u)=l_{u}$ is continuous in its domain.
\end{lem}
\begin{proof}
It can be easily seen that if $l_{u}<\infty$, then
$[0,l_{u}]=\{\mu\in\bbR^{+} | \mu u\in \Omega\}$. In fact, $\langle
l_{u}u,v\rangle =-a$, i.e. if $v=(v_1,\ldots ,v_l )$ then $l_{u}u$
belongs to the hyperplane of $\bbR^{l}$ given by the equation
$v_{1}x_{1}+\ldots+v_{l}x_{l}=-a$. On the other hand $l_{u}u$ also
belongs to the line in $\bbR^{l}$ through the origin in the
direction of $u=(u_1,\ldots ,u_l )$. Therefore $l_u u$ is solution
of the matrix equation $M_{u}X=B$, where
\begin{displaymath}
M_u=\left[\begin{array}{cccc}
v_1 & v_2 &\ldots &v_l\\
u_2-u_1 &\; &\; &\; \\
\; &u_3-u_2 &\; &\; \\
\; &\; &\ddots &\;\\
\; &\; &\; &u_l-u_{l-1}
\end{array}\right],
\indent B=\left[\begin{array}{clcr}-a\\0\\
\vdots \\0\end{array}\right]
\end{displaymath}
and $X$ is the column vector of unknowns.

Let $M_{u}^{i}$ be the matrix obtained exchanging the $i$-th
column of $M_{u}$ by $B$. Using Cramer's Rule, and the fact
$l_{u}=|l_{u} u|$, we obtain the following expression for $f(u)$
\begin{displaymath}
f(u)=l_u= {\left( \sum_{i=1}^{r} {\left(
\frac{detM_{u}^{i}}{detM_{u}}\right)}^{2} \right)}^{\frac{1}{2}}
\end{displaymath}
This proves the result.
\end{proof}

\begin{lem}\label{bounded}
The set $\cP_{\varphi}$ is bounded if and only if there does not
exist $u\in M_{\bbR}$, $u\neq 0$ such that $\bbR^{+}u\subset
\cP_{\varphi}$.
\end{lem}
\begin{proof}
Obviously if $\cP_{\varphi}$ is bounded the conclusion is clear. Let
us prove the other direction. Consider the sets
\begin{displaymath}
\Omega_{i}=\{u\in M_{\bbR} \mid \langle u,v_{i} \rangle\geq
-\varphi_i \}
\end{displaymath}
for all $i\in\cup_{\sigma\in\varDelta}\; [\sigma]=\{1,\ldots
,r\}$. For each $\sigma\in\varDelta$ we have $\gamma_{\sigma
,\varphi}=\cap_{i\in [\sigma ]}\;\Omega_{i}$. Then
$\cP_{\varphi}=\cap \Omega_{i}$ where the intersection is taken
over all the $i\in\{1,\ldots ,r\}$. Let
\begin{equation}\label{lui}
l_{u}^{i}=sup\{\mu\in\bbR^{+} \mid \mu u\in \Omega_{i}\}
\end{equation}
for $u\in \bbR^{l}$ and $i\in\{1,\ldots ,r\}$. Then for each
$i\in\{1,\ldots ,r\}$ we define maps $f_{i}:S\longrightarrow
\bbR^{+}$ where $f_{i}(u)=l_{u}^{i}$ for every $u\in S$. Each
$f_{i}$ is continuous in its domain by Lemma ~\ref{continuous}.
Now we define the map $F:S\longrightarrow\bbR^{+}$ where
\begin{equation}\label{F}
F(u)=\mbox{inf}\{ l_{u}^{i} \mid i\in \{1,\ldots ,r\}\}
\end{equation}
for all $u\in S$. By \cite{Mu}, Chapter 2, Section 18, Exercise 8,
$F$ is continuous in its domain. By hypothesis,  there is no $u\in
M_{\bbR}$, $u\neq 0$ such that $\bbR^{+} u\subseteq
\cP_{\varphi}$. Therefore the  domain of $F$ is $S$. By \cite{Co},
Chp. II 5.8, $F(S^{r-1} )$ is compact in $\bbR^{+}$ and by the
Heine-Borel theorem it is also bounded. Therefore we can find a
positive integer $N$ such that $F(u)\leq N$ for all $u\in S$. If
$y\in \cP_{\varphi}$, then $y=\mu u$ for some $u\in S$ and by
~\eqref{F} and ~\eqref{lui}, $\mu\leq F(u)\leq N$. Therefore
$\cP_{\varphi}$ in included in the sphere of radius $N$.
\end{proof}

\begin{prop}\label{half-space}
 The set $\cP_{\varphi}$ is a polytope
 if and only if the fan $\varDelta $ is not
contained in a half-space.
\end{prop}
\begin{proof}
Suppose there exists $u\in M_{\bbR}$, $u\neq 0$ such that the fan
$\varDelta$ is contained in the half-space $H_{u}$. Then given
$i\in\{1,\ldots ,r\}$ and for any $\mu\in\bbR^{+}$ it holds that
$\langle\mu u,v_{i}\rangle\geq 0\geq -\varphi_i$. Thus
$\bbR^{+}u\subseteq \cP_{\varphi}$ and hence $\cP_{\varphi}$ is
not bounded.

Suppose now that $\cP_{\varphi}$ is not bounded. Then by Proposition
~\ref{bounded} there exists $u\in M_{\bbR}$, $u\neq 0$ such that
$\bbR^{+}u\subseteq \cP_{\varphi}$. By ~\eqref{gamma-sigma-m}, we
have $\mu \langle u,v_{i}\rangle\;\geq -\varphi_i$ for every
$i\in\{1,\ldots ,r\}$ and every $\mu\in\bbR ^{+}$. Then $\langle
u,v_{i}\rangle\geq\; 0$ for all $i\in\{1,\ldots ,r\}$. Therefore the
fan $\varDelta$ is contained in the half-space $H_{u}$.
\end{proof}

%%%%%%%%%%%%%%%%%%%%%%%%%%%%%%%%%%%%%%%%%%%%%%%%%%%%%%%%%%%%%%%%%%%%
\begin{rem}\label{f}
A fan $\varDelta$ is not contained in a half-space if and only if
the intersection of the dual cones is zero,
$\cap_{\sigma\in\varDelta} \sigma^{\vee} =0$. Moreover, given a
set of generators $\{v_1,\ldots ,v_r\}$ of cones in $\varDelta
(1)$ and $0\neq u\in M_{\bbR}$, the following statements are
equivalent:
\begin{enumerate}
\item $\varDelta$ is contained in the half-space $H_u$.

\item The vector $u$ is a solution in $\bbR^l$ of the system of
linear inequalities $\{ \langle x, v_i\rangle \geq 0, i=1,\ldots
r\}$.

\item The set $\bbR_{\geq 0} \{v_1,\ldots ,v_r\}$ is a convex cone
contained in the half-space $H_u$.
\end{enumerate}
\end{rem}

There are well known algorithms to find a nonzero solution of a
system of linear inequalities $\{ \langle x, v_i\rangle \geq 0,
i=1,\ldots , r\}$. Because we need to use it in section
\ref{special}, we give next an algorithm to find such a solution
in the language of half-spaces. The algorithm is based on the
following lemma.

\begin{lem}\label{addvector}
Let $\sigma$ be a cone in $N_{\bbR}$ and let $v\in N_{\bbR}$. Then
$\cC = \bbR_{\geq 0} (\sigma\cup \{v\})$ is a convex cone if and
only if $\cC\subseteq H_u$ for some $u\in M_{\bbR}$ such that
$\langle u, v\rangle =0$.
\end{lem}
\begin{proof}
Let us suppose that $ \{v_1,\ldots ,v_r\}$ is a set of generators
of $\sigma$. Then  $ \{v_1,\ldots ,v_r, v\}$ is a set of
generators of $\cC$. If $\cC$ is a cone, by \cite{Sch}, Theorem
7.1 it is the intersection of the half-spaces $H_u$, $u\in
M_{\bbR}$ containing $ \{v_1,\ldots ,v_r, v\}$ and such that $\{
v\in N_{\bbR} \mid \langle u,v\rangle = 0\}$ is spanned by $l-1$
linearly independent vectors from $ \{v_1,\ldots ,v_r, v\}$. Let
$H_{u_0}$ be one the this half-spaces such that $\langle
u_0,v\rangle = 0$. This proves the result.
\end{proof}

\begin{al}\label{hsp1}
{\sc hsp1}($\{w_{1},\ldots ,w_t\}$ ,$\sigma$ , $u$) {\sf Given} a
convex cone $\sigma$ contained in the half-space $H_u$ with $0\neq
u\in M_{\bbR}$ and a set of vectors $\{w_{1},\ldots ,w_t\}$ in
$N_{\bbR}$, the algorithm {\sf decides }  whether $\cC=\bbR_{\geq
0}(\sigma\cup \{w_{1},\ldots ,w_t\})$ is contained in a half-space
and if the answer is affirmative it {\sf returns } a nonzero
vector $u_0$ such that $\cC$ is contained in $H_{u_0}$.
\begin{enumerate}
\item $i:=1$, $\cC:=\sigma$.

\item While $i\leq t$ do
\begin{itemize}
\item[2.1] If $w_i\in\cC$ then i:=i+1

\item[2.2] else take $u\in M_{\bbR}$ such that $\langle u,
w_i\rangle =0$. If $\cC\subseteq H_{\pm u}$ then $\cC:=\bbR_{\geq
0} (\cC\cup \{w_i\})$, $i:=i+1$, $u_0:=\pm u$ else return {\sf
''Not contained in a half-space"}.
\end{itemize}
\item Return $u_0$ and {\sf ''Contained in a half-space"}.
\end{enumerate}
\end{al}

Let $\cV=\{v_1,\ldots ,v_r\}$ be a set of generators of a fan
$\varDelta$ in $N$. If $\cV$ has rank less than $l$ then
$\varDelta$ is contained in a half-space. Otherwise $\cV$ contains
a subset $\cW$ with $l$ linearly independent vectors. By
\cite{Sch}, Theorem 7.1, $\bbR_{\geq 0} \cW$ is a cone contained
in the half-space $H_u$ for some $u\in M_{\bbR}$ such that
$\cW\subseteq H_u$ and  $\{ v\in N_{\bbR} \mid \langle u,v\rangle
= 0\}$ is spanned by $l-1$ vectors of $\cW$. Algorithm ~\ref{hsp1}
applied to ($\cV\setminus\cW $, $\bbR_{\geq 0} \cW$, $u$) decides
whether $\varDelta $ is contained in a half-space.

%%%%%%%%%%%%%%%%%%%%%%%%%%%%%%%%%%%%%%%%%%%%%%%%%%%%%%%%%%%%%%%%%%%%
\section{Multigradings}\label{multigrading}
%%%%%%%%%%%%%%%%%%%%%%%%%%%%%%%%%%%%%%%%%%%%%%%%%%%%%%%%%%%%%%%%%%%%

The Laurent polynomial ring $S$ is the ring of regular functions on
$X$, $\cO (X)$. We consider the action of $H$ on $\cO (T)$ (or $\cO
(X)$) given by right translation, see \cite{MR}, (12). This
convention implies that $x_i$ has weight $\eta_i$.

We say that a Laurent polynomial ring $S$ is {\em multigraded} by
$A$ when it has been endowed with a degree map $deg
:\bbZ^n\rightarrow A$. Multigradings of polynomial rings are
defined in \cite{MS}, \S 8.1 in the same manner.

Using the matrix $L$ that gives the action of $H$ on $X$, we
define a degree map as follows. The restriction map
$\bbX(T)=Hom(T,k^{\times})\longrightarrow \bbX (H)$ is given by
left multiplication by $L$. If we identify $\bbX(T)$ with $\bbZ^n$
we have an exact sequence (see lemma 4.1)
\begin{equation}\label{multigrading_sequence}
0\longleftarrow A\longleftarrow \bbZ^n\longleftarrow
\FRAK{K}\longleftarrow 0
\end{equation}
which determines a multigrading of S by $A$. Given $\lambda
=(\lambda_1,\ldots ,\lambda_n)$, we write
$x^{\lambda}=x_1^{\lambda_1}\cdot\ldots\cdot x_n^{\lambda_n}$. The
middle semigroup homomorphism $deg :\bbZ^n\rightarrow A$ gives a
multigrading of $S$ by $A$ sending each monomial $x^\lambda$ to
its {\em degree} $deg(\lambda )= L\lambda $. We will say that the
{\em multigrading of $S$ by $A$ is given by the matrix $L$}. Let
$Q=deg(\bbN^r\times\bbZ^s)$ denote the semigroup of $A$ generated
by $deg(x_1),\ldots ,deg(x_r), \pm deg(x_{r+1}),\ldots ,\pm
deg(x_n)$.

%%%%%%%%%%%%%%%%%%%%%%%%%%%%%%%%%%%%%%%%%%%%%%%%%%%%%%%%%%%%%%%%%%%%

The next lemma proves that ~\eqref{multigrading_sequence} is an
exact sequence. Let $l=n-p$.

\begin{lem}\label{deg}
The degree map $deg$ is surjective and there exists an $n\times l$
matrix $K$ whose columns are a $\bbZ$-basis of $\FRAK{K}$.
\end{lem}
\begin{proof}
By ~\cite{A}, Theorem 12.4.3, there exist matrices $U\in GL_m(\bbZ)$
and $V\in GL_n(\bbZ)$ such that
\begin{equation}\label{smith}
L'=ULV=\left[\begin{array}{ccccccc}
f_1   &0  &\ldots &0     &0&\ldots&0\\
0     &f_2&\;     &0     &\;& \;  &\;\\
\vdots&\; &\ddots&\vdots &\vdots&\;&\vdots\\
0  &\ldots&\;     &f_m   &0&\ldots &0
\end{array}\right],
\end{equation}
with $f_i\neq 0$ for all $i=1,\ldots ,m$. Consider the $n\times l$
matrix
\begin{equation}
E=\left[\begin{array}{cc}
0&0\\
0&D\\
I_{l'}&0
\end{array}\right]
\end{equation}
with $D$ as in section \ref{action}, the integer $l'=l-t$ and
$I_{l'}$ the $l'\times l'$ identity matrix.

The columns of $E$ are a $\bbZ$-basis of the kernel of the group
homomorphism  $\bbZ^n\longrightarrow A$ given by right
multiplication by $L'$. Let $K:=VE$, then $LK=0$. Therefore the
columns of $K$ are a $\bbZ$-basis of the kernel of $L$. We can write
$\FRAK{K}$ as $K\bbZ^l$. Then $K$ is a presentation of
$deg(\bbZ^n)$. By \cite{A}, Proposition 5.12 we also have a
presentation $E=V^{-1}K$ of $deg(\bbZ^n)$. We conclude that
$deg(\bbZ^n)=A$.
\end{proof}
%%%%%%%%%%%%%%%%%%%%%%%%%%%%%%%%%%%%%%%%%%%%%%%%%%%%%%%%%%%%%%%%

We are ready to introduce the graded components of $S$ that we are
computing in this paper. For $a\in A$, let $S_a$ denote the vector
space of homogeneous polynomials having degree $a$ in the
$A$-grading, that is

\begin{equation}
S_a=span\{x^\lambda\in S\mid L\lambda =a\}.
\end{equation}
Therefore $S=\bigoplus_{a\in Q} S_a$ since $S_a$ is empty if
$a\notin Q$ and $S_a\cdot S_b = S_{a+b}$, $a,b\in Q$. Note that
the subring $S^H$ of $S$ of invariants under the action of $H$
equals the semigroup ring
\begin{equation}\label{ceroComp}
S_0=k[\FRAK{K}\cap (\bbN^r\times\bbZ^s)].
\end{equation}

%%%%%%%%%%%%%%%%%%%%%%%%%%%%%%%%%%%%%%%%%%%%%%%%%%%%%%%%%%%%%%%%%%%
\section{Construction of fans compatible with the multigrading}\label{compFan}
%%%%%%%%%%%%%%%%%%%%%%%%%%%%%%%%%%%%%%%%%%%%%%%%%%%%%%%%%%%%%%%%%%%

The concept of fan compatible with the multigrading of a
polynomial ring was given in \cite{MS}, section 10.3. Next, we
extent the definition so that it applies to multigradings of
Laurent polynomial rings. In this section, we give an algorithm to
construct a fan compatible with a multigrading from the matrix $L$
describing the action of $H$ on $X$.

Applying the contravariant functor $Hom(\_, \bbZ)$ to the sequence
~\eqref{multigrading_sequence} we obtain the exact sequence
\begin{equation}\label{dual}
0\longrightarrow Hom(A,\bbZ)\longrightarrow \bbZ^n\longrightarrow
\FRAK{K}^{\vee}\longrightarrow Ext^1(A,\bbZ)\longrightarrow 0.
\end{equation}

By the previous lemma ~\ref{deg}, $\FRAK{K}$ is a free abelian
group of dimension l. Then we identify $\FRAK{K}$ and
$\FRAK{K}^{\vee}= Hom(\FRAK{K},\bbZ)$ with $\bbZ^l$. We fix the
columns of the matrix $K$ given by lemma ~\ref{deg} as a
$\bbZ$-basis of $\FRAK{K}$.

The middle morphism $\varphi : \bbZ^n\longrightarrow
\FRAK{K}^{\vee}\equiv \bbZ^l$ of the sequence ~\eqref{dual} is
given by right multiplication by the matrix $K$. Let $e_i$ be the
$i$-th standard basis vector of $\bbZ^n$ and call $v_i =\varphi
(e_i)$ the $i$-th row vector of $K$, $1\leq i\leq n$. In
particular, $\varphi$ is onto when $H$ is an algebraic torus, see
\cite{R}, paragraph following Lemma 2.2. We call $\{v_1,\ldots,
v_r\}$ a {\em set of vectors associated to the action} of $H$ on
$X$.

We call a fan $\varDelta $ in $N$ {\em compatible} with the
multigrading of $S$ by $A$ if its set of one dimensional cones
equals $\varDelta (1) =\{\tau_{v_i} \mid i=1,\ldots ,r\}$. The cones
of $\varDelta$ are subcones of $\cC=\bbR_{\geq 0}\{v_1,\ldots,
v_r\}$.

%%%%%%%%%%%%%%%%%%%%%%%%%%%%%%%%%%%%%%%%%%%%%%%%%%%%%%%%%%%%%%%%%%%

The proof of lemma ~\ref{deg} provides an algorithm to compute a set
of vectors associated to the action of $H$ on $X$.

\begin{al}\label{compatible} {\sc associated-vectors}($L$,$D$)
{\sf Given} the $m\times n$ matrix $L$ associated to the action of
$H$ on $X=k^r\times (k^{\times})^s$ and the $t\times t$ matrix $D$
given in section \ref{action}, the algorithm {\sf returns} a set
of vectors $\{v_1,\ldots ,v_r\}$ associated to the action of $H$
on $X$.
\begin{enumerate}
\item Compute the Smith normal form of $L$; that is, compute
invertible matrices $U$ and $V$ such that $U\cdot L\cdot V$ is the
concatenation of a diagonal matrix and a zero matrix.

\item Let $l':=n-m$, and
\begin{equation} E:=\left[\begin{array}{cc}
0&0\\
0&D\\
I_{l'}&0
\end{array}\right].
\end{equation}

\item $K:=V\cdot E$.

\item Return the set $\{v_1,\ldots ,v_r\}$ of the first $r$ row
vectors of $K$.
\end{enumerate}
\end{al}

%%%%%%%%%%%%%%%%%%%%%%%%%%%%%%%%%%%%%%%%%%%%%%%%%%%%%%%%%%%%%%%%%

%%%%%%%%%%%%%%%%%%%%%%%%%%%%%%%%%%%%%%%%%%%%%%%%%%%%%%%%%%%%%%%%%
\section{Special form of the matrix $L$}\label{specialmatrix}
%%%%%%%%%%%%%%%%%%%%%%%%%%%%%%%%%%%%%%%%%%%%%%%%%%%%%%%%%%%%%%%%%

For some of the constructions in the next sections we need the
matrix $L$ (introduced in section \ref{action}) to have the
special form presented in the next lemma. Now suppose $H$ acts on
$X$ via the matrix $L$ and set
\begin{equation}
\Sigma_L=\{\lambda \in\bbN^r\times\bbZ^s \mid deg(\lambda )=0_A\},
\end{equation}
where $0_A$ is the class of zero in $A$. We show that in some
cases we can obtain a matrix $L'$ equivalent to $L$ with the
special form and giving an action of $H$ on $X$ with the same set
of invariants.

Given a matrix $M$ we will write $\textsl{subm}(M,r_1\ldots r_2,
c_1\ldots c_2)$ to denote the submatrix of $M$ obtained selecting
rows $r_1$ through $r_2$ and columns $c_1$ through $c_2$ of $M$.

\begin{lem}
If $\eta_{r+1},\ldots ,\eta _n$ are linearly independent, there
exist matrices $\Gamma\in GL_m (\bbZ)$, $\Delta\in GL_n (\bbZ)$ such
that
\begin{enumerate}
\item  $L'=\Gamma\cdot L\cdot\Delta$ has the block matrix form
\begin{equation}\label{newL}
\left[\begin{array}{cc}L_1&dI_p\\L_3&L_4\end{array}\right],
\end{equation}
where $I_p$ is the $p\times p$ identity matrix and $d$ is a nonzero
positive integer.

\item $\Sigma_L$ is isomorphic to $\Sigma_L'$.
\end{enumerate}
\end{lem}
\begin{proof}
\begin{enumerate}
\item Since $L$ has rank $m=p+t$ there exists a matrix $\Delta \in
GL_n(\bbZ )$ such that the last $p$ columns of $L\cdot\Delta$ are
linearly independent and such that $\Delta x\in
\bbN^r\times\bbZ^s$ for any $x\in\bbN^r\times\bbZ^s$. To be more
precise, if $p\leq s$ then $\Delta$ is the identity matrix since
$\eta_{r+1},\ldots ,\eta _n$ are linearly independent whereas if
$p>s$ then $\Delta$ is a permutation of the first $r$ columns of
$L$.

Furthermore, we can find a matrix $\Gamma_1\in GL_m(Z)$ permuting
the rows of $L\cdot\Delta$ such that $L_2= subm(\Gamma_1\cdot
L\cdot\Delta, 1\ldots p, l+1\ldots n)$ has rank $p$. Let $d=det(
L_2)$. Then $\Gamma_1=d L_2^{-1}$ is matrix in $GL_p (\bbZ )$. Let
\begin{displaymath}
\Gamma = \left[\begin{array}{cc}\Gamma_1&0\\0&D\end{array}\right]
\end{displaymath}
then $\Gamma\cdot L\cdot \Delta$ has the desired form
~\eqref{newL}.

\item Since $\Delta$ permutes only the first $r$ columns of $L$
and $\Gamma 0_A=0_A$ we can define the following isomorphism
\begin{equation}
\Sigma_L\longrightarrow \Sigma_{L'}
\end{equation}
given by $x\mapsto \Delta^{-1}x$.
\end{enumerate}
\end{proof}

From the proof of statement 1. in the previous lemma we derive the
following algorithm.

\begin{al} {\sc special-matrix}($L$,$D$)
{\sf Given} the $m\times n$ matrix $L$ associated to the action of
$H$ on $X=k^r\times (k^{\times})^s$, where the columns
$\eta_{r+1},\ldots ,\eta _n$ are linearly independent, and the
$t\times t$ matrix $D$ given in section \ref{action}, the
algorithm {\sf returns} an $m\times n$ matrix $L'$ with the
special from ~\eqref{newL}.
\begin{enumerate}
\item $p:=m-t$, $l:=n-p$.

\item If $p> s$ then permute the first $r$ columns of $L$ to
obtain $L_c$ so that\\ $rank(subm(L_c,1\ldots n,l+1\ldots
r))=p-s$.

\item $L:=L_c$. Permute the rows of $L$ to obtain $L_r$ so that
$rank(subm(L_r,1\ldots p, l+1\cdots n))= p$.

\item $L:=L_r$. $L_2:=subm(L, 1\ldots p, l+1\ldots n)$, $d:=det(
L_2)$, $\Gamma_2:=dL_2^{-1}$.

\item $\Gamma :=
\left[\begin{array}{cc}\Gamma_2&0\\0&D\end{array}\right]$.

\item Return $L':=\Gamma\cdot L$.
\end{enumerate}
\end{al}

%%%%%%%%%%%%%%%%%%%%%%%%%%%%%%%%%%%%%%%%%%%%%%%%%%%%%%%%%%%%%%%%%

%%%%%%%%%%%%%%%%%%%%%%%%%%%%%%%%%%%%%%%%%%%%%%%%%%%%%%%%%%%%%%%%%
\section{Determining the special case $S_0=k$}\label{special}
%%%%%%%%%%%%%%%%%%%%%%%%%%%%%%%%%%%%%%%%%%%%%%%%%%%%%%%%%%%%%%%%%

To describe the graded components of $S$ we need to decide in the
first place whether they are finite or infinite dimensional. For
that matter, it is enough to decide whether $S_0$ is finite
dimensional as we will see in section \ref{polytopes}. In this
section, we give an algorithm to decide whether $S_0 =k$ or
equivalently whether $S_0$ is finite dimensional.

We define next an isomorphism allowing us to go from points in the
lattice $M\equiv\bbZ^l$ to points in the lattice $\FRAK
{K}\subset\bbZ^n$ and we use it to prove the next proposition. There
is a natural bilinear pairing $(\;\;,\;):\bbX(T)\times
\bbY(T)\longrightarrow \bbZ$ defined by the requirement that
$(a\circ b)(\lambda )=\lambda^{(a,b)}$ for all $a\in \bbX(T) $,
$b\in \bbY(T)= Hom(k^{\times }, T)$ and $\lambda\in k^{\times }$. We
have an isomorphism $\omega :M\rightarrow \FRAK {K}$ given by
\begin{equation}\label{w}
\langle u,\varphi (b)\rangle =(\omega (u),b)
\end{equation}
for all $u\in M$, $b\in \FRAK {K}$.

\begin{thm}\label{fanS0}
The following statements are equivalent.
\begin{enumerate}
\item The only polynomials of degree 0 are the constants; i.e.
$S_0=k$.

\item The $k$-vector space $S_0$ is finite dimensional.

\item Any fan $\varDelta$ compatible with the multigrading of $S$
given by $A$ is not contained in a half-space.
\end{enumerate}
\end{thm}
\begin{proof}
The proof is analogous to \cite{R}, Lemma 4.2. Let $\varDelta$ be
a fan compatible with the multigrading of $S$ given by $A$. Let
\begin{align}\label{phi}
\phi _{\sigma } &=\{\lambda\in \FRAK{K} \mid (\lambda ,e_{i})\geq 0
\mbox{ for all }i \in [\sigma ]\}.
\end{align}
Then by ~\eqref{ceroComp},
$S_0=k[\cap_{\sigma\in\varDelta}\phi_{\sigma }]$. Therefore $S_0$
is finite dimensional if and only if
$\cap_{\sigma\in\varDelta}\phi_{\sigma }$ is a finite set.
Furthermore $\omega (\cap_{\sigma\in\varDelta} \Lambda_{\sigma
})=\cap_{\sigma\in\varDelta}\phi_{\sigma }$ and
$\cap_{\sigma\in\varDelta} \Lambda_{\sigma }$ is a finite set if
and only if $\cap_{\sigma\in\varDelta} \sigma^{\vee}$ is bounded.
This happens if and only if $\cap_{\sigma\in\varDelta}
\sigma^{\vee}=0$ which is equivalent to statements 1. and 3..
\end{proof}

An $A$-grading of $S$ verifying any of the equivalent conditions
given in the previous proposition is called {\em positive}.

\begin{prop}
The following are necessary conditions for the grading of $S$ by $A$
to be positive.
\begin{enumerate}
\item $p>s$.

\item $\eta_{r+1},\ldots ,\eta_n$ are linearly independent.
\end{enumerate}
\end{prop}
\begin{proof}
\begin{enumerate}
\item Let $\{v_1, \ldots ,v_r\}$ be a set of vectors associated to
the action of $H$ on $X$. If $p\leq s$ or equivalently $l\geq r$
then $\cC =\bbR_{\geq 0}\{v_1, \ldots ,v_r\}$ is a cone and
equivalently $\varDelta$ is contained in a half-space. By theorem
~\ref{fanS0} the result follows.

\item Let $\hat{G}$ be the m dimensional algebraic torus, then
$H\subset\hat{G}$. Consider the action of $\hat{G}$ on $X$ given by
the matrix $L$, then $S^{\hat{G}}\subseteq S^H$. If
$\eta_{r+1},\ldots ,\eta_n$ are linearly dependent, ~\cite{MR},
Lemma 4.1(1) implies that $S^{\hat{G}}$ is not equal to $k$ and
therefore $S_0$ is not finite dimensional.
\end{enumerate}
\end{proof}

Let us suppose that $\eta_{r+1},\ldots ,\eta_n$ are linearly
independent. Then we can assume that $L$ has the special form
~\eqref{newL}. Let $\cV=\{v_1,\ldots ,v_r\}$ be the output of
algorithm ~\ref{compatible} applied to $L_1$ and $D$, and let
$\varDelta$ be a fan compatible with the multigrading of $S$ by
$A$ given by $L$.

For $1\leq i\leq l$, let $\rho_i$ be the restriction of $\eta_i$ to
the subtorus $G$ of $H$. These characters can be thought of as the
columns of $L_1$.

\begin{lem}
If $S_0=k$ then $\rho_i\neq 0$ for all $i=1,\ldots ,l$.
\end{lem}
\begin{proof}
By \cite{R}, Lemma 4.2 if $\rho_i=0$ for some $i=1,\ldots ,l$ then
$\varDelta$ is contained in a half-space and by proposition
~\ref{fanS0} the result follows.
\end{proof}

The converse of the previous lemma does not hold but we can modify
the matrix $L$ to get an action of $H$ on $X$ whose only
invariants are the constant polynomials whenever $\rho_i\neq 0$
for all $i=1,\ldots ,l$. Given a set $I\subseteq \{1,\ldots ,r\}$
define
\begin{equation}
v_i^I=\left\{\begin{array}{c}
-v_i\mbox{ if }i\in I\\
v_i\mbox{ if }i\notin I
\end{array}\right., i=1,\ldots ,r.
\end{equation}
Let $\varDelta_I$ be a fan in $N$ with $\varDelta_I
(1)=\{\tau_{v_i^I} |i=1,\ldots ,r\}$. For $1\leq i\leq n$, set
\begin{equation}\label{LI}
\varsigma_i=\left\{ \begin{array}{c}
-\eta_i\mbox{ if }i\in I\\
\eta_i\mbox{ if }i\notin I
\end{array}\right. .
\end{equation}
Let $L_I$ be the matrix with columns $\varsigma_1,\ldots
,\varsigma_n$. Then $H_I$ denotes the group $H$ acting on $X$ by the
matrix $L_I$. Then the fan $\varDelta_I$ is compatible with the
multigrading of $S$ by $A$ given by $L_I$.

We explain next how to obtain a set $I$ so that $S^{H_I}=k$
whenever $\rho_i\neq 0$ for all $i=1,\ldots ,l$. Let us call
$b_{i,j}$ the entries of $L_1$, $i=1,\ldots ,p$, $j=1,\ldots ,l$.

\begin{lem}\label{lin}
When the matrix $L$ is of the special kind ~\eqref{newL}, then
$v_1,\ldots ,v_l$ are linearly independent.
\end{lem}
\begin{proof}
By Lemma ~\ref{deg}, $LK=0$. Then for $j=l+1,\ldots ,n$
\begin{equation}
dv_j=-\sum_{i=1}^l b_{j-l,i}v_i.
\end{equation}
Thus $v_{l+1},\ldots ,v_n$ belong to the $\bbR$-span of
$v_1,\ldots ,v_l$. By lemma ~\ref{deg} $K$ has rank $l$ and the
result follows.
\end{proof}

By Lemma ~\ref{lin}, $\cB=\{v_1,\ldots ,v_l\}$ is a basis of
$N_{\bbR}$ and  respect to $\cB$ the vectors $v_{l+1},\ldots ,v_n$
have coordinates
\begin{equation}\label{coordinates}
v_j=-\frac{1}{d}(b_{j-l,1},\ldots ,b_{j-l,l}),\indent j=l+1,\ldots
,n.
\end{equation}
Let $v_1^*,\ldots ,v_l^*$ be the dual basis of $\cB$. Then
\begin{equation}\label{coordinates2}
\langle v_i^*, v_j\rangle =\frac{-1}{d}b_{j-l,i}
\end{equation}
for $i=1,\ldots ,l$, $j=l+1,\ldots ,r$.

Let us suppose that $p>s$ or equivalently $l<r$. Given
$j\in\{l+1,\ldots ,r\}$, let
\begin{equation}\label{Iplus}
I_j^+=\{i\in\{1,\ldots ,l\} \mid \frac{-1}{d}b_{j-l,i} \;>0\},
\end{equation}
\begin{equation}\label{Iminus}
I_j^-=\{i\in\{1,\ldots ,l\} \mid \frac{-1}{d}b_{j-l,i} \;<0\}.
\end{equation}

\begin{lem}
If $\rho_i\neq 0$ for all $i=1,\ldots ,l$ then there exists
$J\in\{l+1,\ldots ,r\}$ such that
\begin{equation}\label{J}
\cup_{j=1}^J (I_j^+\cup I_j^-)=\{1,\ldots ,l\}.
\end{equation}
\end{lem}
\begin{proof}
Otherwise there exists $i\in \{1,\ldots ,l\}$ such that
\begin{equation*}
\frac{-1}{d}b_{j-l,i}=0 ,\; j=l+1,\ldots ,r,
\end{equation*}
then $\rho_i=0$.
\end{proof}

This lemma ensures that the next algorithm terminates.

\begin{al}\label{I} {\sc positivity-set}($L_1$, $d$)
{\sf Given} the matrix $L_1$, where the columns $\rho_i\neq 0$ for
all $i=1,\ldots ,l$, and the positive integer $d$ in the special
form of the matrix $L$ giving the action of $H$ on $X$, the
algorithm {\sf returns} a chain of subsets
$\cI_{l+1}\subseteq\ldots\subseteq\cI_J$ of $\{1,\ldots ,l\}$ such
that the grading of $S$ by $A$ given by $L_{\cI_J}$ is positive.
\begin{enumerate}
\item $J:=l+1$.

\item Compute $I_J^-$ and $I_J^+$ using ~\eqref{Iplus} and ~\eqref{Iminus}.

\item $\cI_J:=I_J^+$.

\item While $\cup_{j=l+1}^J (I_j^-\cup I_j^+)\neq \{ 1,\ldots
,l\}$ then
\begin{enumerate}
\item[4.1] $I_J^0:=\{ 1,\ldots ,l\}\setminus I_J$.

\item[4.2] Compute $I_{J+1}^-$ and $I_{J+1}^+$.

\item[4.3] $\cI_{J+1}:=\cI_J\cup((\cap_{j=l+1}^{J}I_j^0)\cap
I_{J+1}^+)$.

\item[4.4] $J:=J+1$.
\end{enumerate}
\item Return $\cI_{l+1},\ldots ,\cI_J$.
\end{enumerate}
\end{al}

The next proposition shows that the grading of $S$ by $A$ given by
$L_{\cI_J}$, with $\cI_J$ as in the output of algorithm ~\ref{I}, is
positive.

\begin{prop}\label{not in half-space}
Let us suppose that $\rho_i\neq 0$ for all $i=1\ldots l$ and let
$\cI_J$ be the subset of $\{1,\ldots ,l\}$ obtained by {\sc
positivity-set}($L_1$,$d$), then the grading of $S$ by $A$ given
by $L_{\cI_J}$ is positive.
\end{prop}
\begin{proof}
It can be proved as in \cite{R}, Proposition 4.5 that
$\varDelta_\cI $ is not contained in a half-space. The result
follows by theorem ~\ref{fanS0}.
\end{proof}

If $\cI_J=\emptyset$ then the previous proposition ensures that
$\varDelta = \varDelta_{\emptyset}$ is not contained in a
half-space. On the other hand, if $\cI_J\neq\emptyset$ the fan
$\varDelta$ might still be contained in a half-space. We give next
an algorithm to check whether $\varDelta$ is contained in a
half-space using the output set of algorithm ~\ref{I}. The
algorithm is based on the following fact.

\begin{lem}
Let $\{\cI_{l+1} ,\ldots ,\cI_J\}$ be the output of {\sc
positivity-set}($L_1$, $d$). Given $k\in \{1,\ldots ,J-(l+1)\}$,
if $\cI_{J-k}=\emptyset$ and $\cI_{J-k+1}\neq \emptyset$ then
$\cC= \bbR_ {\leq 0}\{v_1,\ldots ,v_{J-k+1}\}$ is contained in
$H_{v_i^*}$ for all $i\in \cI_{J-k+1}$.
\end{lem}
\begin{proof}
Given $i\in \cI_{J-j+1}$ then $i\in (\cap_{j=l+1}^{J-k}I_j^0)\cap
I_{J-k+1}^+$. Thus by ~\eqref{coordinates2}, ~\eqref{Iplus},
~\eqref{Iminus} and 4.1 in algorithm ~\ref{I} we have $\langle
v_i^*, v_j\rangle\geq 0$ for all $j=l+1,\ldots ,J-k+1$. This
proves the result.
\end{proof}

\begin{al}\label{hsp2}
{\sc hsp2}($\{\cI_{l+1} ,\ldots ,\cI_J\}$, $\cV$) {\sf Given}
$\{\cI_{l+1} ,\ldots ,\cI_J\}$ the output set of algorithm
~\ref{I} and $\cV=\{v_1,\ldots ,v_r\}$ a set of vectors associated
to the action of $H$ on $X$, the algorithm {\sf decides} whether
$\cC=\bbR_{\geq 0} \cV$ is contained in a half-space and if the
answer is affirmative it {\sf returns} a vector $0\neq u\in
M_{\bbR}$ such that $\cC$ is contained in $H_u$.
\begin{enumerate}
\item $\cC:=\bbR_{\geq 0}\{v_1,\ldots ,v_l\}$.

\item If $\cI_{l+1},\ldots ,\cI_J$ are nonempty sets then
\begin{itemize}
\item[2.1] $\cC:=\bbR_{\geq 0}(\cC\cup\{v_{l+1}\})$, $u:=v_i^*$ with $i\in
\cI_{l+1}$.

\item[2.2] If $J=l+1$ return $u$ and {\sf '' Contained in a half-space"}
else return {\sc hsp1}($\{v_{l+2},\ldots ,v_J\}$, $\cC$ , u).
\end{itemize}

\item If $\cI_J=\emptyset$ then return  {\sf '' Not contained in a half-space"}.

\item Let $1\leq i\leq J-(l+1)$ be the smallest such $\cI_{J-i}=\emptyset$ and $\cI_{J-i+1}\neq\emptyset$.
\begin{itemize}
\item[4.1] $\cC:=\bbR_{\geq 0}(\cC\cup\{v_{l+1},\ldots
,v_{J-i+1}\})$ and $u:=v_i^*$ with $i\in \cI_{J-i+1}$.

\item[4.2] If $J-i+1=J$ return $u$ and {\sf '' Contained in a half-space"}.

\item[4.3] Return {\sc hsp1}($\{v_{J-i+2},\ldots ,v_J\}$, $\cC$ ,
u).
\end{itemize}

\end{enumerate}
\end{al}

From the previous results we derive the following algorithm to
test whether an $A$-grading is positive.

\begin{al} {\sc positivity-test}($L$)
{\sf Given} the $m\times n$ matrix $L$ associated to the action of
$H=G\times \bbF$ on $X=k^r\times (k^{\times})^s$ where $G$ is a
torus of dimension $p$, the algorithm {\sf decides} whether the
$A$-grading of $S$ by $L$ is positive and if the answer is
negative it {\sf returns} $\cI\subset\{1,\ldots l\}$ such that
$A$-grading of $S$ by $L_{\cI}$ is positive.
\begin{enumerate}
\item If $p>s$ or $\eta_{r+1},\ldots ,\eta_n$ are linearly
dependent return {\sf ''The grading is not positive"}.

\item
$\left[\begin{array}{cc}L_1&dI_p\\L_3&L_4\end{array}\right]:=${\sc
special-matrix}($L$,$D$).

\item If any of the columns of $L_1$ is $0$ then return {\sf ''The
grading is not positive"}.

\item $\{\cI_{l+1},\ldots ,\cI_J\}:=${\sc
positivity-set}($L_1$,$d$).

\item If {\sc hsp2}($\cI_{l+1} ,\ldots ,\cI_J$) returns {\sf ''Not
contained in a half-space"} then return {\sf''The grading is
positive"} else return $\cI_J$ and {\sf ''The grading is not
positive"}.
\end{enumerate}
\end{al}

%%%%%%%%%%%%%%%%%%%%%%%%%%%%%%%%%%%%%%%%%%%%%%%%%%%%%%%%%%%%%%%%%%%%%
\section{Polyhedral description of graded components of S}\label{polytopes}
%%%%%%%%%%%%%%%%%%%%%%%%%%%%%%%%%%%%%%%%%%%%%%%%%%%%%%%%%%%%%%%%%%%%%

In this section we describe the graded components of $S$ in terms
of polyhedra. We distinguish to main cases depending on the
dimension.

Let $\varDelta$ be a fan compatible with the multigrading of $S$
by $A$ given by $L$. Given $a\in Q$, there exists
$\varphi=(\varphi_1,\ldots ,\varphi_n)\in\bbN^r\times\bbZ^s$ such
that $L\varphi=a$. Consider the polyhedron $\cP_{\varphi}$ defined
by ~\eqref{polytope}.

\begin{lem}\label{Sa}
Let $a\in Q$ then for any $\varphi=(\varphi_1,\ldots
,\varphi_n)\in\bbN^r\times\bbZ^s$ such that $L\varphi=a$ it holds
\begin{equation}\label{Sa2}
S_a=x^{\varphi}k[\omega (\cP_{\varphi}\cap\bbZ^l)],
\end{equation}
with $\omega$ as in ~\eqref{w}.
\end{lem}
\begin{proof}
For each $a\in Q$ we have $S_a=k[(\varphi + \FRAK{K})\cap
(\bbN^r\times\bbZ^s)]$. Given $\sigma\in\varDelta$ let
\begin{align}\label{phim}
\phi _{\sigma ,a} &=\{\varphi +\mu\in \varphi +\FRAK{K} \mid (\mu
,e_{i})\geq -\varphi_i\mbox{ for all } i\in [\sigma ]\}.
\end{align}
Then $S_a=k[\displaystyle{ \cap_{\sigma\in\varDelta}} \phi
_{\sigma ,a}]$. Define the sets
\begin{equation}
\psi_{\sigma ,a}=\{\lambda\in \FRAK{K} \mid (\lambda ,e_{i})\geq
-\varphi_i, \mbox{ for all }i\in [\sigma ]\}.
\end{equation}
Thus $\phi_{\sigma ,a}=\varphi +\psi_{\sigma ,a}$ and
$k[\phi_{\sigma ,a} ]= x^{\varphi } k [\psi_{\sigma ,a} ]$.
Therefore $S_a=x^{\varphi} k[ \cap_{\sigma} \psi_{\sigma ,a} ]$.
Let
\begin{equation}
\Lambda_{\sigma ,a}=\{x\in M\mid \langle x,v_{i}\rangle\geq
-\varphi_i\;\mbox{for all}\; i\in [\sigma ]\}.
\end{equation}
Then $\psi_{\sigma ,a}=\omega (\Lambda_{\sigma ,a})$ and $\omega
(\cap_{\sigma\in\varDelta} \Lambda_{\sigma ,a })= \cap_{\sigma}
\psi_{\sigma ,a}$. Finally $\cap_{\sigma\in\varDelta}
\Lambda_{\sigma ,a }=\cP_{\varphi}\cap \bbZ^l$.
\end{proof}

The following theorem gives two more properties of positive
gradings.

\begin{thm}
The following statements are equivalent:
\begin{enumerate}
\item There exits $a\in A$ such that the $k$-vector space $S_a$ is
finite dimensional.

\item For all $a\in A$, the $k$-vector space $S_a$ is finite
dimensional.

\item Any fan compatible with the multigrading of $S$ given by $A$
is not contained in a half-space.
\end{enumerate}
\end{thm}
\begin{proof}
Let $\varDelta$ be a fan compatible with the multigrading of $S$
by $A$. Given $\varphi\in\bbN^r\times\bbZ^s$, by proposition
~\ref{half-space} the set $\cP_{\varphi}\cap \bbZ^l$ is finite if
and only if $\varDelta$ is not contained in a half-space. Given
$a\in Q$ and $\varphi\in\bbN^r\times\bbZ^s$ such that $L\varphi
=a$ by the previous lemma $S_a$ is finite dimensional if and only
if $\cP_{\varphi}\cap \bbZ^l$ is finite. This proves the result.
\end{proof}

%%%%%%%%%%%%%%%%%%%%%%%%%%%%%%%%%%%%%%%%%%%%%%%%%%%%%%%%%%%%%%%%%
\subsection{Finite dimensional case}

At this point we can determine the graded components in the finite
dimensional case.

\begin{thm}\label{thmFin}
Let us suppose that the grading of $S$ by $A$ given by $L$ is
positive. Given $a\in Q$, then the dimension of $S_a$ equals the
number of lattice points inside the polytope $\cP_{\varphi}$ for any
$\varphi$ such that $L\varphi =a$.
\end{thm}
\begin{proof}
By proposition  ~\ref{half-space}, the polyhedron $\cP_{\varphi }$
is bounded. By lemma ~\ref{Sa} the graded component $S_a$ is
spanned over $k$ by the finite set of monomials $x^u$ such that
\begin{displaymath}
u\in \varphi+\omega (\cP_{\varphi}\cap \bbZ^l).
\end{displaymath}
This proves the result
\end{proof}

\begin{rem}
{\sc polyhedral description} The proof of theorem ~\ref{thmFin}
provides an algorithm to determine $S_a$ using a polytope.
\end{rem}

Let us suppose that the grading of $S$ by $A$ given by $L$ is
positive, that is $S_0=k$. Then by \cite{MR}, Lemma 4.1(1) the
vectors $\eta_{r+1},\ldots ,\eta_n$ are linearly independent and
we can assume that $L$ is of the special kind ~\eqref{newL}. Let
$\{v_1,\ldots ,v_r\}$ be a set of vectors associated to the action
of $H$ on $X$. By lemma ~\ref{lin}, the set
$\bar{\cB}=\{-v_1\ldots ,-v_l\}$ is an $N_{\bbR}$ basis. Let us
consider the $r\times l$ matrix $P$ whose $i$-th row is the row
vector of the coordinates of $v_i$ in the $N_\bbR$ basis
$\bar{\cB}$, $i=1\ldots r$. In the dual basis of $\bar{\cB}$ the
polytope $\cP_{\varphi}$ equals the set
\begin{equation}\label{poli}
\{x\in \bbR^l| Px\leq \varphi_i, i=1,\ldots ,r\}.
\end{equation}

%%%%%%%%%%%%%%%%%%%%%%%%%%%%%%%%%%%%%%%%%%%%%%%%%%%%%%%%%%%%%%%%%
\subsection{Infinite dimensional case}

We describe next the graded components of $S$ for the infinite
dimensional case. Let us suppose that $S_0$ is not finite
dimensional. Let $\cV =\{v_1,\ldots ,v_r\}$ be a set of vectors
associated to the action of $H$ on $X$ and let $\varDelta$ be a
fan compatible with the multigrading of $S$ by $A$. By theorem
~\ref{fanS0}, $\varDelta$ is contained in a half-space and by
remark ~\ref{f}, the set $\cC =\bbR_{\geq 0} \{v_1,\ldots ,v_r\}$
is a cone. Let $\cC^{\vee}$ be the dual cone of $\cC$. By
\cite{Sch}, Theorem 16.7 there exists a Hilbert basis of
$\cC^{\vee}$.
\begin{prop}
Let $S_0$ be infinite dimensional and let $\cH=\{w_1,\ldots
,w_h\}$ be a Hilbert basis of $\cC^{\vee}$. Then
\begin{displaymath}
S_0 = k[x^{\omega (w_1)},\ldots ,x^{\omega (w_h)}].
\end{displaymath}
\end{prop}
\begin{proof}
 For each
$\sigma\in\varDelta$, let $\phi _{\sigma }$ be defined as in
~\eqref{phi} and $\Lambda_{\sigma}$ as in ~eq\ref{lambda}. Then
$\cap_{\sigma\in\varDelta} \phi_{\sigma} = \FRAK{K}\cap
(\bbN^r\times\bbZ^s)$,  $\cap_{\sigma\in\varDelta}
\Delta_{\sigma}= \cC^{\vee}\cap M$ and $\omega
(\cap_{\sigma\in\varDelta} \Delta_{\sigma}) =
\cap_{\sigma\in\varDelta} \phi_{\sigma}$. Thus $\{ \omega
(w_1),\ldots ,\omega (w_h)\}$ is a minimal generating set of the
partially ordered set $\FRAK{K}\cap (\bbN^r\times\bbZ^s)$ and this
proves the result.
\end{proof}

Every graded component of $S$ is infinite dimensional. Given $a\in
Q$ let $\varphi=(\varphi_1,\ldots
,\varphi_n)\in\bbN^r\times\bbZ^s$ such that $L\varphi=a$. The next
result will be used to determine $S_a$.

By proposition ~\ref{half-space}, the polyhedron $\cP_{\varphi }$
is not bounded. By \cite{Sch}, Corollary 7.1b then $\cP_{\varphi
}=\cP +\cC^{\vee}$ for some polytope $\cP$.

\begin{lem}\label{descomposition}
Let $\cH=\{w_1,\ldots ,w_h\}$ be a Hilbert basis of $\cC^{\vee}$,
let $B =\{ \sum_{i=1}^h \alpha_i w_i \mid 0\leq \alpha_i \leq 1\}$
and let $\cP$ be a polytope such that $\cP_{\varphi } = \cP
+\cC^{\vee}$. Then

\begin{equation}
\cP_{\varphi }\cap\bbZ^l = ((\cP
+B)\cap\bbZ^l)+(\cC^{\vee}\cap\bbZ^l)
\end{equation}

\end{lem}
\begin{proof}
We prove the nontrivial inclusion $\cP_{\varphi }\cap\bbZ^l
\subseteq ((\cP +B)\cap\bbZ^l)+(\cC^{\vee}\cap\bbZ^l)$. Given
$p\in \cP_{\varphi }\cap\bbZ^l$ then $p=q+c$ with $q\in Q$ and
$c\in \cC^{\vee}$. We can write $c=\sum \alpha_i w_i$ with
$\alpha_i\geq 0$. Consider the integral vector
 $c'=\sum \lfloor \alpha_i\rfloor w_i$, thus $b=c-c'\in B$. We
 can write $p=q+b+c'$ with $q+b=p-c'\in(\cP +B)\cap\bbZ^l$ and $c'\in \cC^{\vee}\cap\bbZ^l$.
\end{proof}

Given a set of generators of a cone, we can use the algorithm given
in \cite{H}, section 5.5 to compute a Hilbert basis of the dual
cone.

\begin{thm}\label{thmInf}
For all $a\in Q$, $S_a$ is finitely generated as an $S_0$-module.
\end{thm}
\begin{proof}
Let us consider the cone $\bar{\cC}$ generated by vectors
$(v_i,\varphi_i)\in\bbZ^{l+1}$, $i=1,\ldots ,r$ and
$(0_{\bbZ^l},1)\in \bbZ^{l+1}$ with $0_{\bbZ^l}$ the zero vector
in $\bbZ^l$. By \cite{Sch}, Theorem 16.7 there exists a Hilbert
basis $\cH$ of the dual cone $\bar{\cC}^{\vee}$ of $\bar{\cC}$
which equals
\begin{displaymath}
\{(x,\lambda)\in \bbR^l\times\bbR_{\geq 0}\mid \langle
x,-v_i\rangle -\lambda\varphi_i\leq 0\}
\end{displaymath}
The set $\cH_0=\{(x,\lambda)\in\cH\mid \lambda=0\}$ is a Hilbert
basis of $\cC^{\vee}$. Let $\cP$ be the convex hull of the
$x\in\bbZ^l$ with $(x,1)\in \cH$. Now $x\in \cP_{\varphi}$ if and
only if $(x,1)\in \bar{\cC}^{\vee}$. Then $\cP_{\varphi } = \cP
+\cC^{\vee}$.

By lemma ~\ref{Sa} and lemma  ~\ref{descomposition} we conclude
that $S_a$ is generated as an $S_0$-module by the elements $x^u$
such that $u\in \varphi+\omega ((\cP +B)\cap\bbZ^l)$.
\end{proof}

\begin{rem}
{\sc polyhedral description} The proof of theorem ~\ref{thmInf}
provides an algorithm to determine $S_a$ using polyhedra.
\end{rem}

\end{document}